\documentclass{article}
\usepackage{graphicx}
\usepackage{amsmath}
\usepackage{amssymb}
\usepackage{amsthm}

\title{On Cusp Solutions to a Prescribed Mean Curvature Equation}
\author{Alexandra K. Echart \& Kirk E. Lancaster
                                  \\
                       Department of Mathematics, Statistics \& Physics \\
                            Wichita State University \\
                            Wichita, Kansas, 67260-0033}

\def\Real{{\rm I\hspace{-0.2em}R}}

\newcommand\myeq{\mathrel{\overset{\makebox[0pt]{\mbox{\normalfont\tiny\sffamily def}}}{=}}}

\date{ }

\begin{document}
\maketitle

\vspace*{3mm}

\begin{abstract}
The nonexistence of ``cusp solutions'' of prescribed mean curvature boundary value problems in $\Omega\times\Real$  
when $\Omega$  is a domain in $\Real^{2}$ is proven in certain cases and an application to radial limits at a corner is mentioned.  
\end{abstract}
\newtheorem{thm}{Theorem}
\newtheorem{prop}{Proposition}
\newtheorem{cor}{Corollary}
\newtheorem{lem}{Lemma}
\newtheorem{rem}{Remark}
\newtheorem{example}{Example}

\section{Introduction}
Let $\Omega$  be a domain  in ${\Real}^{2}$  with locally Lipschitz boundary  and ${\cal O}=(0,0)\in \partial\Omega$ 
and $H\in C^{1,\beta}\left(\overline{\Omega}\times \Real\right),$  for some $\beta\in (0,1).$     
Let polar coordinates relative to ${\cal O}$  be denoted by $r$ and $\theta$  and let $B_{\delta}({\cal O})$ 
be the open ball in $\Real^{2}$  of radius $\delta$ about ${\cal O}.$
We shall assume there exists a $\delta^{*}>0$  and $\alpha \in \left(0,\pi\right)$  such that 
$\partial \Omega \cap B_{\delta^{*}}({\cal O})$   consists of two smooth arcs  ${\partial}^{+}\Omega^{*}$  and 
$\partial^{-}\Omega^{*}$, whose tangent lines approach the lines $L^{+}: \:  \theta = \alpha$  
and  $L^{-}: \: \theta = - \alpha$, respectively, as the point ${\cal O}$ is approached and 
for each $\theta\in (-\alpha, \alpha),$  there exists an $r(\theta)>0$  such that 
$\{(r\cos(\theta),r\sin(\theta)): 0<r<r(\theta)\} \subset \Omega.$  
Set $\Omega^{*} = \Omega \cap B_{\delta^{*}}({\cal O}).$

\begin{figure}[ht]
\centering
\includegraphics{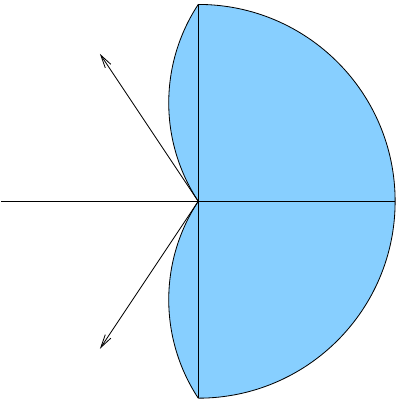}
\caption{The domain $\Omega^{*}$  \label{ZERO}}
\end{figure}

Consider a solution $f\in C^{2}(\Omega)$  of the prescribed mean curvature equation
\begin{equation}
\label{PMC}
{\rm div}(Tf)(x,y)  =   2H(x,y,f(x,y)) \ \ \ \ {\rm for} \ \  (x,y)\in\Omega^{*}, 
\end{equation}
which satisfies the conditions  
\begin{equation}
\label{Bounds}
\sup_{(x,y)\in\Omega^{*}} |f(x,y)| < \infty \ \ \ \  {\rm and} \ \ \ \  \sup_{(x,y)\in\Omega^{*}} |H(x,y,f(x,y))| < \infty,  
\end{equation}
where $Tf = \frac{\nabla f}{\sqrt{1 + |\nabla f|^{2}}};$  
examples of such functions might arise as solutions of a Dirichlet or contact angle boundary value problem for (\ref{PMC}).
We are interested in the radial limits of $f:$   
\begin{equation}
\label{Rads}
Rf(\theta) \myeq \lim_{r\downarrow 0} f(r\cos(\theta),r\sin(\theta)), \ \ \ -\alpha<\theta<\alpha. 
\end{equation} 
When $\lim_{\partial^{+}\Omega^{*}\ni (x,y)\to {\cal O} } f\left(x,y\right)$  exists,  we define 
$Rf(\alpha)$  to be this limit and when $\lim_{\partial^{-}\Omega^{*}\ni (x,y)\to {\cal O} } f\left(x,y\right)$  
exists,  we define $Rf(-\alpha)$  to be this limit.   

There are examples in which the radial limits do not exist for any $\theta\in (-\alpha,\alpha)$  (\cite{Lan:89, LS1}).  
For solutions of boundary value problems which satisfy appropriate conditions, $Rf(\theta)$  can be proven to exist 
for $\theta\in [-\alpha,\alpha]\setminus J,$  where $J$  is a countable subset of $(-\alpha,\alpha)$  
(e.g. \cite{NoraKirk1, NoraKirk2, Lan:88, Lan:91, Lan:12, LS1, LS2}). 
We know of no examples in which $J\neq\emptyset$  and we ask if $J=\emptyset$  always holds; 
this is related to the existence of {\it cusp solutions}.

A {\it cusp solution} for (\ref{PMC}) is a domain $\Lambda\subset\Real^{2}$  and a solution $f$  of (\ref{PMC}) in 
$\Lambda$  such that $\partial\Lambda\setminus \{ {\cal O},A,B\}=\Gamma_{1}\cup\Gamma_{2}\cup\Gamma_{3},$  where 
$A,B,{\cal O}$  are distinct points on $\partial\Lambda,$   $\Gamma_{1},$  $\Gamma_{2}$  and $\Gamma_{3}$  are 
disjoint, smooth (open) arcs with respective endpoints $\{A,{\cal O}\},$  $\{B,{\cal O}\}$  and $\{A,B\},$
$\Gamma_{1}$  and $\Gamma_{2}$  are tangent at ${\cal O}$  (so  $\overline{\Lambda}$  has an 
``outward'' cusp at ${\cal O};$  see Figure \ref{ONE}, which has a cusp at $(0,0)$), 
$f(x,y)=c_{j}$  when $(x,y)\in \Gamma_{j}$  ($j=1,2$),  $c_{1} < c_{2},$  
and, for each $c\in (c_{1},c_{2}),$  the level curves $\{ (x,y)\in\Lambda \ : \ f(x,y)=c\}$  are tangent at ${\cal O}$
(e.g. \S 5 of \cite{LS1}).  
(Capillary surfaces in cusp regions were studied in \cite{AS,Scholz:03}.)  
In cases where cusp solutions do not exist, we know that $J=\emptyset.$

In \cite{LS1, LS2}, the nonexistence of cusp solutions is proven when 
(a) $H\in C^{1,\delta}\left(\overline{\Omega}\times\Real\right)$  ($\delta\in (0,1)$)  and $H(x,y,z)$  is 
strictly increasing in $z$  for each $(x,y)\in\overline{\Omega}$  or (b) $H$  is real-analytic.  
The proof in \cite{LS1} for case (a) involves a ``local'' argument while that for case (b) involves a ``global'' argument 
which shows that (\ref{Bounds}) is violated.  
Using a ``local'' argument, we shall prove 

\begin{thm}
\label{Theorem1}
Suppose $\Omega$  is a domain  in ${\Real}^{2}$  with locally Lipschitz boundary,  ${\cal O}=(0,0)\in \partial\Omega$  and   
$H\in C^{1,\beta}\left(\overline{\Omega^{*}}\times \Real\right)$  for some $\beta\in (0,1).$   
Let $f\in C^{2}(\Omega^{*})$  satisfy (\ref{PMC}) and (\ref{Bounds}).  
Suppose $H(x,y,z)$  is weakly increasing in $z$  for $(x,y)$  in a neighborhood of $(0,0).$  
Then $f$  cannot have a cusp solution (i.e. there is no ``cusp region'' $\Lambda\subset\Omega$  such that 
$\left(\Lambda,f\right)$  is a cusp solution).
\end{thm} 

\noindent We can exclude cusp solutions when  $H$  vanishes  in the ``cusp direction,'' which we may assume is the direction of the positive 
$x-$axis  (see Figure \ref{ONE}). 

\begin{thm}
\label{Theorem2}
Suppose $\Lambda$  is a cusp domain  in ${\Real}^{2},$  $\partial\Lambda$  is tangent to $\vec i$  at ${\cal O},$    
$H\in C^{1,\beta}\left(\overline{\Lambda}\times \Real\right)$    for some $\beta\in (0,1),$  
$f\in C^{2}(\Lambda)$  satisfies (\ref{PMC}) and (\ref{Bounds}) and there exists a $\delta>0$  such that 
\[
H(x,0,z)=0 \ \ \ {\rm for} \ \ \  
(x,z)\in [0,\delta]\times [\liminf_{\Lambda\ni (x,y)\to {\cal O}} f(x,y), \limsup_{\Lambda\ni (x,y)\to {\cal O}} f(x,y)].
\]
Then $(\Lambda,f)$  cannot be a cusp solution.  
\end{thm} 

\noindent What can we say when $H(x,y,z)$  is strictly decreasing in $z$?  
Unfortunately, as the following example illustrates, we cannot exclude cusp solutions in this case, even when $H$  is real-analytic; 
a ``global'' argument (like \cite{LS1}, page 176) is required to exclude cusp solutions when $H$  is real-analytic.  
Thus, for example, the reasoning in 3B of \cite{AS} cannot be used when $\kappa<0.$  

\begin{example}
Consider the cone ${\cal C}=\{X(\theta,t) : 0\le \theta\le \frac{\pi}{2}, 0< t<\infty\},$  where
\[
X(\theta,t) = t(\cos(\theta), \sin(\theta)-1, 1).
\]
Set $\Lambda=\{ t(\cos(\theta), \sin(\theta)-1) : 0< \theta< \frac{\pi}{2}, 1< t< 2\}$  and    
${\cal S}={\cal C}\cap \left(\Real^{2}\times [1,2]\right).$
A straightforward computation shows that the mean curvature (with respect to the upward normal) is 
\[
H(\theta,t)=\frac{3-2\sin(\theta)}{2t\left(1+(1-\sin(\theta))^{2}\right)^{3/2}};
\]
in other words, $H(x,y,z)=\frac{z^2-2yz}{2\left(y^{2}+z^{2}\right)^{3/2}}.$   Now $\frac{y}{z}=\sin(\theta)-1 \in [-1,0]$  
and $x=0$  iff $\theta=\pi/2;$  another calculation yields  
\[
2\frac{\partial H}{\partial z}(x,y,z)=-\frac{z^{3}}{\left(y^{2}+z^{2}\right)^{5/2}}
\left(1-4\left(\frac{y}{z}\right)-2\left(\frac{y}{z}\right)^{2}+2\left(\frac{y}{z}\right)^{3}\right)
< 0.
\]
Finally observe that ${\cal S}$  is the graph of a cusp solution and satisfies (\ref{Bounds}) in $\Lambda.$ 
\end{example}

The hypotheses of \cite{NoraKirk1} include the assumption that $H$  satisfies one of the conditions which 
guarantees that cusp solutions  do not exist; the following Corollary is a consequence of Theorem \ref{Theorem1} and \cite{NoraKirk1}.    
(A second corollary, similar to Corollary \ref{Corollary}, follows by applying Theorem \ref{Theorem1} to 
Theorems 1 \& 2 of \cite{NoraKirk2}.)

\begin{cor} (\cite{NoraKirk1})
\label{Corollary}
Suppose $\Omega,$  $f$  and $H$  satisfy the hypotheses of Theorem \ref{Theorem1} and either   
\begin{itemize}
\item[(i)] $\alpha \in \left(\frac{\pi}{2},\pi\right)$  or  
\item[(ii)] $\alpha \in \left(0,\frac{\pi}{2}\right]$  and one of $Rf(\alpha)$ or $Rf(-\alpha)$  exists.
\end{itemize}
Then $Rf(\theta)$  exists for each $\theta\in (-\alpha,\alpha)$  and $Rf\in C^{0}\left( (-\alpha, \alpha) \right).$  
If $Rf(\alpha)$  exists, then $Rf\in C^{0}\left( (-\alpha, \alpha] \right).$
If $Rf(-\alpha)$  exists, then $Rf\in C^{0}\left( [-\alpha, \alpha) \right).$
\end{cor}

\section{Proof of Theorem \ref{Theorem1}}

\begin{figure}[ht]
\centering
\includegraphics{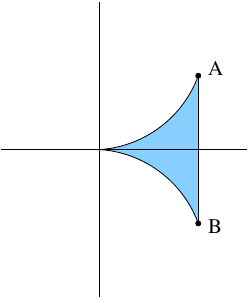}
\caption{The cusp domain $\Lambda$  \label{ONE}}
\end{figure}

Suppose $\left(\Lambda,f\right)$  is a cusp solution and $\Lambda\subset \{(x,y)\in \Real^{2}: 0<x<a,\ |y|<x\},$   
$c_{1}<c_{2}$  and the $c-$level curves of $f$  in $\Lambda$  are tangent to the positive $x-$axis at ${\cal O}$  
for $c_{1}\le c \le c_{2},$  for some $a>0$  (see Figure \ref{ONE}).  
Since $H\in C^{1,\beta}\left(\overline{\Omega}\times \Real\right),$  $f\in C^{3}(\Omega)$  and, 
as in \cite{LS1,LS2}, there exist an (open) rectangle $R_{0}=(0,a)\times (c_{1},c_{2})$  
and $g\in C^{3}\left(R\right)$   (where $R=\overline{R_{0}}$)   
such that the graph of $f$  over $\Lambda,$  $\cal{G},$  is the set $\{(x,g(x,z),z) : (x,z)\in R_{0}\}$  
(i.e. $z=f(x,y)$  iff $y=g(x,z)$  for $(x,z)\in R_{0}$  and $(x,y)\in \Lambda$)  
and $g(0,z)=\frac{\partial g}{\partial x}(0,z)=0$  for $c_{1}\le z\le c_{2}.$   
We may assume that $|\nabla g(x,z)|\le 1$  for $(x,z)\in R.$

The (upward) unit normal to the graph of $f,$  $\cal{G},$   is 
\[
\vec N(x,y,z)=\frac{(-f_{x}(x,y),-f_{y}(x,y),1)}{\sqrt{1+f_{x}^{2}(x,y)+f_{y}^{2}(x,y)}}
\]
and ${\rm div}(Tf)(x,y) = 2\vec H(x,y,z) \cdot \vec N(x,y,z)$  for $(x,y,z)\in \cal{G},$  
where $2\vec H$  is the mean curvature vector of $\cal{G}.$  
Then 
\[
{\rm sgn}(g_{z}(x,z)) \vec N(x,y,z)=\frac{(g_{x}(x,z),-1,g_{z}(x,z)))}{\sqrt{1+g_{x}^{2}(x,z)+g_{z}^{2}(x,z)}}. 
\]
Since ${\rm div}(Tg)=2\vec H \cdot (-g_{x},1,-g_{z})/\sqrt{1+g_{x}^{2}+g_{z}^{2}},$  we see that 
\[
{\rm div}(Tg)(x,z) = 2\vec H(x,y,z) \cdot (-{\rm sgn}(g_{z}(x,z)))\vec N(x,y,z) \ \ \ \ {\rm for} \ (x,y,z)\in \cal{G}.
\]
(Of course, if $g_{z}(x,z)=0$  for some $(x,z)\in R$  with $x>0,$  then $\cal{G}$  has a horizontal unit normal at 
an interior point of $\Omega,$  which contradicts our hypothesis $f\in C^{2}(\Omega);$  hence $g_{z}(x,z)\neq 0$
when $(x,z)\in R$  with $x>0.$)  

Let us assume ${\rm sgn}(g_{z}(x,z))={\rm sgn}(f_{y}(x,g(x,z)))=+1$  for $(x,z)\in R$  with $x>0;$   the opposite choice 
will lead to the same (eventual) conclusion that cusp solutions do not exist.  
Then 
\[
Mg(x,z)=-2H(x,g(x,z),z),
\]
where $Mg = \nabla \cdot Tg = {\rm div}\left(Tg\right).$ 
Suppose there exist a $\delta_{1}>0$  such that $H(x,y,z)$  is weakly increasing in $z$  for each $(x,y)\in\Lambda$  and 
$z\in [c_{1},c_{2}]$  when $x^{2}+y^{2}\le \delta_{1}^{2}.$   We may assume $a\le\delta_{1}.$

Fix  $\epsilon\in \left(0,\frac{1}{2}(c_{2}-c_{1})\right)$  and set $\tilde c_{1}=c_{1}+\epsilon$  and 
$\tilde c_{2}=c_{2}-\epsilon;$  notice that $\tilde c_{2}>\tilde c_{1}.$  
Set 
\begin{equation}
\label{Harvest}
g_{j}(x,z):=g\left(x,z+\tilde c_{j}\right) \ \ \ \ {\rm for} \  0\le x\le a,\ -\epsilon \le z \le \epsilon, \ \ \ j=1,2,    
\end{equation}
and define $h=g_{1}-g_{2}.$   

If $h(x_{0},z_{0})=0$  for some $(x_{0},z_{0})\in (0,a]\times [-\epsilon,\epsilon],$  then the graph of $f$  fails the vertical 
line test since $\left(x_{0},y_{0},z_{0}+\tilde c_{1}\right)$  and $\left(x_{0},y_{0},z_{0}+\tilde c_{2}\right)$  
are both points on the graph of $f,$   where $y_{0}=g_{1}(x_{0},z_{0})$  ($=g_{2}(x_{0},z_{0})$).  
Thus $h(x,z)\neq 0$  for all $0<x\le a, \ -\epsilon\le z\le \epsilon.$  
Since ${\rm sgn}(g_{z}(x,z))=+1$  when $(x,z)\in (0,a]\times [-\epsilon,\epsilon],$  we see that $h(x,z)<0$  for all 
$(x,z)\in (0,a]\times [-\epsilon,\epsilon].$  
(This is essentially the argument in \cite{LS1}  (at the bottom of page 175) since $h(0,z)>0$  is the only option available there.)

Define 
\[
K(x,y)=2H(x,y,\tilde c_{1}+\epsilon), \ \ \ 0\le x\le a, \ (x,y)\in\Lambda.
\]
and  $d(x,z)=2H(x,g(x,z),\tilde c_{1}+\epsilon) - 2H(x,g(x,z),z).$
Notice that $d(x,z+\tilde c_{1})\ge 0$  and $d(x,z+\tilde c_{2})\le 0$  when $(x,z)\in [0,a]\times[-\epsilon,\epsilon].$    
Now, for each $j=1,2,$  $g_{j}$  is a solution of the Cauchy problem 
\begin{eqnarray*}
\label{Cauchy3} 
Mg_{j}(x,z) & = & -K(x,g_{j}(x,z))+d(x,z+\tilde c_{j}) \ \ \ \ {\rm for} \ \ (x,z)\in [0,a]\times [-\epsilon,\epsilon]  \\
\label{Cauchy4}
g_{j}(0,z) & = & \frac{\partial g_{j}}{\partial x}(0,z)  =  0 \ \ \ \  {\rm for} \ \ z \in [-\epsilon,\epsilon].
\end{eqnarray*} Then, as in \cite{GT}, pp. 263-4, we have
\begin{eqnarray*}
0 & = & Mg_{1}(x,z)-Mg_{2}(x,z)+2H(x,g_{1}(x,z),z+\tilde c_{1})-2H(x,g_{2}(x,z),z+\tilde c_{2}) \\ 
& = & Lh(x,z)-d(x,z+\tilde c_{1})+d(x,z+\tilde c_{2}),
\end{eqnarray*}
where, setting  $D_{1}:= \frac{\partial}{\partial x}$  and $D_{2}:= \frac{\partial}{\partial z},$   
\begin{equation}
\label{Pizza}
Lh= \sum_{i,j=1}^{2} a^{i,j}D_{ij}h + \sum_{i=1}^{2} b^{i}D_{i}h + ch;
\end{equation}
here 
\begin{equation}
\label{Apple1}
a^{i,j}(x,z)=e^{i,j}(Dg_{1}(x,z)) \ \ \  \ {\rm for} \ \  i,j=1,2,  
\end{equation}
with $e^{1,1}(p,q)= (1+q^2)W^{-3},$  $e^{1,2}(p,q)= e^{2,1}(p,q)= -pqW^{-3},$  $e^{2,2}(p,q)= (1+p^2)W^{-3},$  $W=W(p,q)=\sqrt{1+p^2+q^2},$ 
\begin{equation}
\label{Apple2}
b^{1}(x,z)=\sum_{i,j=1}^{2} D_{ij}g_{2}(x,z)\frac{\partial e^{i,j}}{\partial p}(\xi_{1},(g_{1})_{z}(x,z)),
\end{equation}
\begin{equation}
\label{Apple3}
b^{2}(x,z)=\sum_{i,j=1}^{2} D_{ij}g_{2}(x,z)\frac{\partial e^{i,j}}{\partial q}((g_{2})_{x}(x,z),\xi_{2})
\end{equation}
and 
$c(x,z)= \frac{\partial K}{\partial y}(x,\xi)=2\frac{\partial H}{\partial y}(x,\xi,\tilde c_{1}+\epsilon),$   
for some $\xi$  between $g_{1}(x,z)$  and $g_{2}(x,z),$  $\xi_{1}$  between $(g_{1})_{x}(x,z)$  and $(g_{2})_{x}(x,z)$  and
$\xi_{2}$  between $(g_{1})_{z}(x,z)$  and $(g_{2})_{z}(x,z).$  

Notice that $a^{i,j}\in C^{1}\left(R\right)$  for $i,j=1,2,$   $b^{i}\in L^{\infty}(R)$  for $i=1,2$  and $c\in L^{\infty}(R).$
Now $h(0,z)=\frac{\partial h}{\partial x}(0,z)=0$  for $|z|\le \epsilon$  and 
\begin{equation}
\label{Nice}
Lh(x,z)=d(x,z+\tilde c_{1})-d(x,z+\tilde c_{2})\ge 0, \ \ \  (x,z)\in [0,a]\times [-\epsilon,\epsilon]. 
\end{equation}
From (\ref{Nice}) and the Hopf boundary point lemma (e.g. \cite{GT}, Lemma 3.4), we have 
\[
\frac{\partial h}{\partial x}(0,z)<0 \ \ \ {\rm for \ each} \ \ z\in (-\epsilon,\epsilon)
\]
and this contradicts the fact that $h_{x}(0,z)=0$  if $z\in [-\epsilon,\epsilon].$
Thus we have proven Theorem \ref{Theorem1}.

\begin{rem}
The assumption that $H$  is weakly increasing in $z$  is equivalent to one in the (weak) comparison principle 
(e.g. Theorem 10.1 in \cite{GT}; Theorem 5.1 in \cite{FinnBook}), which plays a critical role here.  
\end{rem}

\section{Proof of Theorem \ref{Theorem2}}

Suppose $\left(\Lambda,f\right)$  is a cusp solution and $\Lambda\subset \{(x,y)\in \Real^{2}: 0<x<a,\ |y|<x\},$   
$c_{1}<c_{2}$  and the $c-$level curves of $f$  in $\Lambda$  are tangent to the positive $x-$axis at ${\cal O}$  
for $c_{1}\le c \le c_{2},$  for some $a>0$  (see Figure \ref{ONE}).   
As before, there exist an (open) rectangle $R_{0}=(0,a)\times (c_{1},c_{2})$  
and $g\in C^{3}\left(R\right)$   such that the graph of $f$  over $\Lambda,$  $\cal{G},$  
is the set $\{(x,g(x,z),z) : (x,z)\in R_{0}\}$  and $g(0,z)=\frac{\partial g}{\partial x}(0,z)=0$  for $c_{1}\le z\le c_{2}.$   
We shall assume that $|\nabla g(x,z)|\le 1$  for $(x,z)\in R.$

Let us assume there exist $\delta\in (0,a]$  and $d_{1},d_{2}\in [c_{1},c_{2}]$  with $d_{1}<d_{2}$  such that 
$H(x,0,z)=0$  for $0\le x\le\delta,$  $d_{1}\le z\le d_{2}.$  
Now  $g_{xx}(0,z)=0$  for all $z\in [c_{1},c_{2}]$  (since  $\triangle g(0,z)=Mg(0,z)=-2H(0,0,z)=0$) and 
\[
H(x,g(x,z),z) = H(x,0,z)  + \frac{\partial H}{\partial y}(x,\xi,z)g(x,z)  = \frac{\partial H}{\partial y}(x,\xi,z)g(x,z)
\]
for some $\xi$  between  $0$  and $g(x,z).$  
We may extend $g$  as an even function in $x$  by setting $g(x,z)=g(-x,z)$  for $-a\le x<0,$  $c_{1}\le z\le c_{2},$  
so that $g\in C^{2}(R\cup R^{-}),$  where $R^{-}=\{(-x,z) : (x,z)\in R\}.$  
Then 
\[
0  =  Mg(x,z)+2H(x,g(x,z),z) =  \tilde Lg(x,z)
\]
where  $a^{1,1}(x,z)=\frac{1+g_{z}^{2}(x,z)}{W^{3}},$ 
$a^{1,2}(x,z)=-\frac{g_{x}(x,z)g_{z}(x,z)}{W^{3}},$  
$a^{2,2}(x,z)=\frac{1+g_{x}^{2}(x,z)}{W^{3}},$
$W(x,z)=\sqrt{1+g_{x}^{2}(x,z)+g_{z}^{2}(x,z)},$
$a^{1,2}=a^{2,1},$
$\tilde c(x,z)=2 H_{y}(x,\xi,z)$  and 
\[
\tilde Lu= \sum_{i,j=1}^{2} a^{i,j}D_{ij}u +  \tilde cu.
\]
Since  $|\nabla g(x,z)|\le 1$  for $(x,z)\in R,$  $\tilde L$  is uniformly elliptic in $R.$  
Notice that $a^{i,j}\in C^{1}\left(R\right)$  for $i,j=1,2$  and  $\tilde c\in C^{0}(R).$  
Since  $g\in C^{2}(R\cup R^{-}),$
Theorems $1^{*}$  and $2^{*}$  of \cite{HW:53} imply 
that for each $z\in (d_{1},d_{2}),$  there exist a natural number $n$  and 
real constants $e_{1}$  and $e_{n},$  not both zero, such that 
\[
g_{x}(\rho\cos(\theta),z+\rho\sin(\theta)) = \rho^{n}\left(e_{1}\cos(n\theta)+e_{2}\sin(n\theta)\right) + o(\rho^{n})
\]
and 
\[
g_{z}(\rho\cos(\theta),z+\rho\sin(\theta)) = \rho^{n}\left(e_{2}\cos(n\theta)-e_{1}\sin(n\theta)\right) + o(\rho^{n})
\]
as $\rho\to 0.$  
Since $g_{x}(0,z)=0$  and $g_{z}(0,z)=0$  for $z\in [c_{1},c_{2}],$  we see that 
\[
e_{1}\cos(n\pi/2)+e_{2}\sin(n\pi/2)=0, \ \ \ \  e_{2}\cos(n\pi/2)-e_{1}\sin(n\pi/2)=0
\]
and so $e_{1}=e_{2}=0.$  This contradicts the fact that at least one of $e_{1}$  or $e_{2}$  is non-zero. 
Thus we have proven Theorem \ref{Theorem2}.

\section{Radial Limits}

When radial limits for (\ref{PMC}) exist, they behave in a different manner than do radial limits of, for example, 
Laplace's equation (e.g. \cite{BH:83}).  In particular,  if $f$  is a solution of (\ref{PMC}) and the radial limits 
$Rf(\theta)$  exist for $\theta\in (-\alpha, \alpha),$   then they behave in one of the following ways:
\begin{itemize}
\item[(i)]  $Rf:(-\alpha,\alpha)\to\Real$  is a constant function (i.e. $f$  has a nontangential limit at ${\cal O}$). 
\item[(ii)] There exist $\alpha_{1}$ and $\alpha_{2}$ so that $-\alpha \leq \alpha_{1}
< \alpha_{2} \leq \alpha$ and $Rf$ is constant on $(-\alpha, \alpha_{1}]$ and
$[ \alpha_{2}, \alpha)$ and strictly increasing or strictly decreasing on
$(\alpha_{1}, \alpha_{2})$.  
\item[(iii)] There exist $\alpha_{1}, \alpha_{L}, \alpha_{R}, \alpha_{2}$ so that
$-\alpha \leq \alpha_{1} < \alpha_{L} < \alpha_{R} < \alpha_{2} \leq \alpha,
\alpha_{R}= \alpha_{L} + \pi$, and $Rf$ is constant on $(-\alpha, \alpha_{1}],
[ \alpha_{L}, \alpha_{R}]$, and $[ \alpha_{2}, \alpha)$ and either  strictly increasing
on $(\alpha_{1}, \alpha_{L}]$ and  strictly decreasing on $[ \alpha_{R}, \alpha_{2})$ or
 strictly decreasing on $(\alpha_{1}, \alpha_{L}]$ and  strictly increasing on $[\alpha_{R},\alpha_{2})$.  
\end{itemize}

\end{document}